\documentclass[11pt,reqno]{amsart}
\usepackage{amssymb}
\usepackage{amsmath, amssymb, amsthm}
\usepackage{mathrsfs}

\newtheorem{theorem}{Theorem}[section]

\newtheorem{lemma}{Lemma}[section]
\newtheorem{remark}{Remark}[section]

\numberwithin{equation}{section}

\begin{document}
\title[Magnetohydrodynamic equations]{Blow-up criterion for the  $3$D non-resistive  compressible Magnetohydrodynamic equations}

\author{Shuai Xi}
\address[S. Xi]{School of Mathematical Sciences, Shanghai Jiao Tong University,
Shanghai 200240, P.R.China}
\email{\tt shuai\_xi@sina.com}

\author{Shengguo Zhu}
\address[S. G. Zhu]{The Institute of Mathematical Sciences, The Chinese University of Hong Kong, Shatin, N.T., Hong Kong;  Department of Mathematics, Shanghai Jiao Tong University,
Shanghai 200240,  P.R.China.}
\email{\tt sgzhu@mail1.ims.cuhk.edu.hk}

\begin{abstract}In this paper, we prove a blow-up criterion in terms of   the magnetic field $H$ and the mass density $\rho$ for the strong solutions to  the  $3$D  compressible isentropic   MHD equations with zero magnetic diffusion and initial vacuum.  More precisely,  we  show   that  the upper bounds  of  $(H,\rho)$  control the possible blow-up (see \cite{olga}\cite{zx}\cite{by}) for  strong solutions.
\end{abstract}

\date{Feb. 28, 2014}

\maketitle

\section{Introduction}
Magnetohydrodynamics is that part of the mechanics of continuous media which studies the motion of electrically conducting media in the presence of a magnetic field. The dynamic motion of  fluid and  magnetic field interact strongly on each other, so the hydrodynamic and electrodynamic effects are coupled.
In $3$D space, the compressible isentropic   MHD equations in a domain $\Omega \subset \mathbb{R}^3$  can be written as
\begin{equation}
\label{eq:1.2j}
\begin{cases}
\displaystyle
H_t-\text{rot}(u\times H)=-\text{rot}\Big(\frac{1}{\sigma}\text{rot}H\Big),\\[6pt]
\displaystyle
\text{div}H=0,\\[6pt]
\displaystyle
\rho_t+\text{div}(\rho u)=0,\\[6pt]
\displaystyle
(\rho u)_t+\text{div}(\rho u\otimes u)
  +\nabla P=\text{div}\mathbb{T}+\text{rot}H\times H.
\end{cases}
\end{equation}

In this system,  $x\in \Omega$ is the spatial coordinate; $t\geq 0$ is the time;  $H=(H^1,H^2,H^3)$ is the magnetic field; $0< \sigma\leq \infty$ is the electric conductivity coefficient;  $\rho\geq 0$ is the mass density;
$u=(u^1,u^2,u^3)\in \mathbb{R}^3$ is the velocity of fluids; $P$ is the pressure  satisfying
\begin{equation}
\label{eq:1.3}
P=A \rho^\gamma, \quad A>0,\quad \gamma > 1, 
\end{equation}
where $A$  is a   constant  and $\gamma$ is the adiabatic index;
$\mathbb{T}$ is the viscosity stress tensor:
\begin{equation}
\label{eq:1.4}
\mathbb{T}=2\mu D(u)+\lambda \text{div}u \mathbb{I}_3, \quad D(u)=\frac{\nabla u+(\nabla u)^\top}{2},
\end{equation}
where $D(u)$ is the deformation tensor,  $\mathbb{I}_3$ is the $3\times 3$ unit matrix, $\mu$ is the shear viscosity coefficient, $\lambda+\frac{2}{3}\mu$ is the bulk viscosity coefficient, $\mu$ and $\lambda$ are both real constants,
\begin{equation}
\label{eq:1.5}
  \mu > 0, \quad 3\lambda+2\mu \geq 0,\end{equation}
which ensures the ellipticity of the   Lam$\acute{\text{e}}$ operator (see (\ref{xc})).
Although the electric field $E$ doesn't appear in  system (\ref{eq:1.2j}),
it is indeed induced according to a relation
$$
E=\frac{1}{\sigma}\text{rot}H- u\times H
$$
by moving the conductive flow in the magnetic field.

The MHD system (\ref{eq:1.2j}) describes the macroscopic behavior of electrically conducting compressible (isentropic) fluids
in a magnetic field. It is reasonable to assume that there is no magnetic diffusion (i.e. $\sigma=+\infty$) when the conducting fluid considered is of a very high conductivity,  which occurs frequently in many cosmical and geophysical problems. Then we need  to consider the following system:
\begin{equation}
\label{eq:1.2}
\begin{cases}
\displaystyle
H_t-\text{rot}(u\times H)=0,\\[6pt]
\displaystyle
\text{div}H=0,\\[6pt]
\displaystyle
\rho_t+\text{div}(\rho u)=0,\\[6pt]
\displaystyle
(\rho u)_t+\text{div}(\rho u\otimes u)
  +\nabla P=\text{div}\mathbb{T}+\text{rot}H\times H,
\end{cases}
\end{equation}
which  is the so called viscous and non-resistive MHD equations (see \cite{dd1}\cite{dd2}\cite{dd3}\cite{dd4}\cite{dd5}).

The aim of this paper is to give a blow-up criterion of strong solutions  to the initial boundary value problem (IBVP):  system (\ref{eq:1.2}) in  a bounded,  smooth domain $\Omega \subset \mathbb{R}^3$ with the initial boundary value  conditions:
\begin{align}
\label{fan1}(H,\rho, u)|_{t=0}=(H_0(x), \rho_0(x),  u_0(x)),\  x\in \Omega; \quad u|_{\partial {\Omega}}=0.
\end{align}

Throughout this paper, we adopt the following simplified notations for the standard homogeneous and inhomogeneous Sobolev space:
\begin{equation*}\begin{split}
 &\|f\|_{W^{m,r}}=\|f\|_{W^{m,r}(\Omega)},\quad  \|f\|_s=\|f\|_{H^s(\Omega)},\quad |f|_p=\|f\|_{L^p(\Omega)},\\[6pt]
&D^{k,r}=\{f\in L^1_{loc}(\Omega): |\nabla^kf|_{r}<+\infty\},\quad  |f|_{D^{k,r}}=\|f\|_{D^{k,r}(\Omega)}=|\nabla^kf|_r, \quad D^k=D^{k,2}, \\[6pt]
& |f|_{D^k}=\|f\|_{D^k(\Omega)}=|\nabla^kf|_2,\quad D^{1}_0=\{f\in L^6(\Omega): |\nabla f|_{2}<\infty \ \text{and}\ f|_{\partial \Omega}=0\},\\
& |f|_{D^1_0}=\|f\|_{D^1_0(\Omega)}=|\nabla^kf|_2,\quad \|(f,g)\|_X=\|f\|_X+\|g\|_X, \quad  \mathbb{A}: \mathbb{B}=\sum_{ij}a_{ij}b_{ij}.
\end{split}
\end{equation*}
A detailed study of homogeneous Sobolev spaces can be found in \cite{gandi}.

As has been observed in Theorem $5.1$ of  \cite{jishan},  which proved the existence of unique local strong solutions  with initial vacuum to IBVP (\ref{eq:1.2})-(\ref{fan1}),  in order to make sure that the  IBVP (\ref{eq:1.2})-(\ref{fan1})   with initial vacuum is well-posed, the lack of a positive lower bound of the initial mass density $\rho_0$ should be compensated with some initial layer compatibility condition on the initial data $(H_0,\rho_0, u_0,P_0)$:
\begin{theorem}\cite{jishan}\label{th5} Let constant $q \in (3,6]$.
If  $(H_0, \rho_0,  u_0,P_0)$ satisfies
\begin{equation}\label{th78}
\begin{split}
& (H_0, \rho_0, P_0) \in H^1\cap W^{1,q},\  \rho_0\geq 0,\    u_0\in D^1_0\cap D^2,
\end{split}
\end{equation}
and the compatibility condition
\begin{equation}\label{th79}
\begin{split}
Lu_0+\nabla P_0- \text{rot} H_0\times H_0=\sqrt{\rho_0} g_1
\end{split}
\end{equation}
for some $g_1 \in L^2$, $P_0=A\rho^\gamma_0$,  and
\begin{equation}\label{xc}
\begin{split}
Lu_0=-\mu\triangle u_0-(\lambda+\mu)\nabla  \text{div} u_0,
\end{split}
\end{equation}
then there exists a  time $T_*$ and a unique solution $(H,\rho,u,P)$ to IBVP (\ref{eq:1.2})-(\ref{fan1}) satisfying
\begin{equation*}\begin{split}
&(H,\rho,P)\in C([0,T_*];H^1\cap W^{1,q}),\ u\in C([0,T_*];D^1_0\cap D^2)\cap L^{2}([0,T_*];D^{2,q}),\\
& u_t\in  L^2([0,T_*];D^1_0),\ \sqrt{\rho}u_t\in L^\infty([0,T_*];L^2).
\end{split}
\end{equation*}
\end{theorem}

Some analogous existence theorems of the unique local strong solutions to the  compressible Navier-Stokes equations   have been previously established by Cho-Choe-Kim in  \cite{CK3}\cite{CK}\cite{guahu}.  In $3$D space, Huang-Li-Xin obtained the well-posedness of global classical solutions with small energy but possibly large oscillations and vacuum for Cauchy problem in  \cite{HX1} or IBVP in \cite{HX2} to the isentropic flow.
For compressible MHD equations, when $0<\sigma<+\infty$,
 the  global smooth solution near the constant state in one-dimensional space was studied in Kawashima-Okada \cite{kawa}; recently, in $3$D space, the similar result to \cite{HX1} has been obtained in Li-Xu-Zhang \cite{mhd}, and the existence of global weak solutions has also been proved in Hu-Wang \cite{hx}.  However, for $\sigma=+\infty$,  as far as we know, there are few results  on the  global existence of strong solutions with initial vacuum.  The non-global existence in the whole space $\mathbb{R}^3$ has been proved in \cite{olga} for the classical solution   to isentropic MHD equations  as follows:
\begin{theorem}\cite{olga} \label{olga}
Assume that $\gamma\geq \frac{6}{5}$, if the momentum $\int_{\mathbb{R}^3} \rho_0 u_0 \text{d} x\neq 0 $, then there exists no global classical solutions  to (\ref{eq:1.2})-(\ref{fan1}) with conserved mass, momentum and total energy.
\end{theorem}
The similar results also can be seen in Yuan-Zhao \cite{by}.

Then these motivate us to consider  that the local strong solutions to (\ref{eq:1.2})-(\ref{fan1})  may cease to exist globally,  or what is the key point to make sure that the solution obtained in Theorem \ref{th5} could become a  global one?
If the blow-up happens,  we want to know  the mechanism of breakdown and the structure of  singularities?
The similar question has been studied for the incompressible Euler equation  by Beale-Kato-Majda (BKM) in their pioneering work \cite{TBK}, which showed that the $L^\infty$-bound of vorticity $\text{rot} u$ must blow up if we assume that the life span of the corresponding strong solution is finite. Later, Ponce  \cite{pc} rephrased the BKM-criterion in terms of the deformation tensor $D(u)$, and the  same result as \cite{pc} has been proved by Huang-Li-Xin \cite{hup} for  compressible isentropic Navier-Stokes equations, which can be shown: if  $0 < \overline{T} < +\infty$  is the
maximum existence time for strong solution, then
\begin{equation}\label{kaka1}
\lim \sup_{ T \rightarrow \overline{T}} \int_0^T |D( u)|_{ L^\infty(\Omega)}\text{d}t=\infty,
\end{equation}
which have also  been proved  for the $3$D  compressible isentropic MHD equations in Zhu \cite{szz}.
We refer reders to  Liu-Liu \cite{mingtao},    Lu-Du-Yao  \cite{duyi} and   Xu-Zhang \cite{gerui}  for more details.

Moreover, for the strong solutions with initial vacuum to   $3$D compressible isentropic Navier-stokes equations, Sun-Wang-Zhang \cite{zif} proved
$$
\lim \sup_{ T \rightarrow \overline{T}} |\rho|_{ L^\infty([0,T]\times \Omega)}=\infty,
$$
under the assumption  (\ref{eq:1.5})  and $\lambda <7\mu$, which has been improved to $\lambda <\frac{29}{3}\mu$  in Wen-Zhu \cite{jiangzhu}.

Recently,  when $0<\sigma<+\infty$, some interesting results have been obtained for the compressible  isentropic MHD system (\ref{eq:1.2j}):
$$
\lim \sup_{ T \rightarrow \overline{T}}( |\rho|_{L^\infty([0,T]; L^{p_1}(\Omega))}+ |H|_{L^\infty([0,T]; L^{p_2}( \Omega))})=\infty,
$$
for any $1<p_1<\infty$ and $p_2=\infty$ under the assumption (\ref{eq:1.5}) with $\lambda<\frac{29}{3}\mu$ in \cite{yongfu}, or for any $\frac{24}{5}\leq p_2 \leq \infty$ and  $p_1=\infty$  under the assumption (\ref{eq:1.5}) with  $ \lambda<3\mu$ in \cite{yu}.

So it is interesting and important  to ask whether the similar conclusions obtained in  \cite{hup}\cite{zif}\cite{yongfu}\cite{yu} also hold for system (\ref{eq:1.2}) with $\sigma=\infty$.  Compared with \cite{yongfu}\cite{yu}, due to the vanishing of magnetic diffusion, there are two significant difficulties in our proof: the strong coupling between $u$ and $H$ and  the lack of smooth mechanism of $H$.  For example, in order to  deal with  the  nonlinear term: magnetic momentum flux density tensor
$$
\frac{1}{2}|H|^2\mathbb{I}_3-H\otimes H
$$
 in momentum equations $(\ref{eq:1.2})_4$,  we need to control the norm $ |\nabla H|_{2}$,  which is difficult to be  bounded by $|D(u)|_{L^1(0,T;L^\infty)}$ because of the issues mentioned above. These are unlike the estimates for $(|\rho|_\infty, |\nabla \rho|_{2})$ which  can be totally determined by $|\text{div}u|_{L^1(0,T;L^\infty)}$ due to the  scalar hyperbolic structure of the continuity equation $(\ref{eq:1.2j})_1$ in  \cite{hup}\cite{zif}, and the estimate for $|\nabla H|_2$ in \cite{yongfu}\cite{yu} based on the strong parabolic structure of magnetic equations when $0<\sigma<\infty$.  So some new arguments need to be introduced to improve the results obtained above   for system  (\ref{eq:1.2j}).

In the following   theorem, under the physical assumption  (\ref{eq:1.5})  and  $3\lambda < 29 \mu$,  we  show that the $L^\infty$ norms of the magnetic field $H$ and  the mass density $\rho$  control the possible blow-up  (see \cite{olga}\cite{zx}) for  strong solutions, which  means that if a solution of the compressible MHD equations is initially regular and loses its regularity at some later time, then the formation of singularity must be caused  by  losing  the upper  bound of   $H$ or $\rho$ as the critical time approaches.

\begin{theorem} \label{th3}
Let   the viscosity coefficients $(\mu, \lambda)$ satisfy
\begin{equation}\label{ass}
 \mu > 0, \quad 3\lambda+2\mu \geq 0,\quad3\lambda < 29 \mu,\end{equation}
 and   $(H_0,\rho_0,u_0, P_0)$ satisfy (\ref{th78})-(\ref{th79}).
 If  $(H,\rho,u,P)$  is a strong solution to IBVP  (\ref{eq:1.2})- (\ref{fan1}) obtained in Theorem \ref{th5},   and $0< \overline{T} <\infty$ is the maximal time of  its existence, then
\begin{equation}\label{eq:2.91}
\lim \sup_{ T \rightarrow \overline{T}}( |\rho|_{L^\infty([0,T]\times \Omega)}+ |H|_{L^\infty([0,T]\times \Omega)})=\infty.
\end{equation}
\end{theorem}

\begin{remark}\label{rrr2}
We introduce the main ideas of our proof   for Theorem \ref{th3}.

\text{I}) First, if $ |\rho|_{L^\infty([0,T]; L^{\infty}(\Omega))}$ and $ |H|_{L^\infty([0,T]; L^{\infty}(\Omega))}$ are bounded, we can obtain a high integrability of velocity $u$, which can be used to control the nonlinear terms (See Lemmas \ref{abs3}-\ref{sk4}). For this purpose,    an  important observation has been shown  in Lemma \ref{sk4} that
\begin{equation}\label{magnet}\begin{cases}
\quad \mathbb{B}=H_t\otimes H+H\otimes H_t\\[6pt]
=(H^iH^k\partial_k u^j+H^jH^k\partial_k u^i-H^i H^j\partial_k u^k)_{(ij)}-\text{div}\big( (H\otimes H)\otimes u \big),\\[6pt]
\quad \mathbb{C}=H\cdot H_t=H\cdot (H \cdot \nabla u-u \cdot \nabla H-H\text{div}u)\\[6pt]
=\big(H\cdot \nabla u\cdot H-\frac{1}{2}|H|^2\text{div}u\big)-\frac{1}{2}\text{div}(u|H|^2),
\end{cases}
\end{equation}
from which, we successfully avoid the difficulty coming from the   strong coupling between the magnetic field and velocity when the magnetic diffusion vanishes.

\text{II}) Second, the next difficulty is to control the mass density $\rho$ and the magnetic field $H$, which both satisfy hyperbolic equations.
To do this, we need to make sure that  the velocity $u$ is bounded in $L^1([0, T]; D^{1,\infty}(\Omega))$.  On the other hand, in order to prove $u\in L^1([0, T]; D^{1,\infty}(\Omega))$,
we have to obtain some priori bounds for $\nabla \rho$ and $\nabla H$.
Furthermore, the magnetic  term  in the momentum equation will bring
extra difficulty to us.  However, via using the argument from \cite{hoff1} and  the structure of the magnetic equations,  in Lemma \ref{ablem:4-1},  we show that
\begin{equation}\label{mgnetc}
\begin{split}
\Lambda=& \int_{\Omega} \dot{u} \cdot \big[\text{div}\big(H\otimes H-\frac{1}{2}|H|^2I_3\big)_t+\text{div}\big(\text{div}\big(H\otimes H-\frac{1}{2}|H|^2I_3\big)\otimes u\big)\big]\text{d}x\\
=& \int_{\Omega}  \partial_{k}u^k H^i H^j  \partial_{j}  \dot{u}^i  \text{d}x+\int_{\Omega} \Big( -\frac{1}{2}\partial_ju^j |H^k|^2 \partial_{i}\dot{u}^i\Big) \text{d}x.
\end{split}
\end{equation}
Then we get the cancelation to the derivatives $(\nabla \rho, \nabla H)$ during our computation, which brings us the desired result.
\end{remark}

The rest of  this paper is organized as follows. In Section $2$, we give some important lemmas which will be used frequently in our proof.    In Section $3$,  we give the proof for the blow-up criterion (\ref{eq:2.91}).

\section{Preliminary}

In this section, we give some important lemmas which will be used frequently in our proof.
The first one is some Sobolev inequalities:

\begin{lemma}\label{gag}
For  $l\in (3,\infty)$, there exists some generic constant $C> 0$ that may depend  on $l$ such that for $f\in D^1_0(\Omega)$,  $g\in D^1_0\cap D^2(\Omega)$  and $h \in W^{1,l}(\Omega)$, we have
\begin{equation}\label{gaga1}
\begin{split}
|f|_6\leq C|f|_{D^1_0},\qquad |g|_{\infty}\leq C|g|_{D^1_0\cap D^2}, \qquad |h|_{\infty}\leq C\|h\|_{W^{1,l}}.
\end{split}
\end{equation}
\end{lemma}

Next we consider the following boundary value problem for the Lam$\acute{\text{e}}$ operator $L$:
\begin{equation}\label{tvd}
\begin{cases}
-\mu \triangle U-(\mu+\lambda)\nabla \text{div} U=F,\quad \text{in}\ \Omega,\\[8pt]
U(t,x)=0, \quad \text{on}\quad \partial \Omega,
\end{cases}
\end{equation}
where $U = (U^1, U^2, U^3), \ F = (F^1, F^2, F^3)$. It is well known that under the assumption $(\ref{eq:1.5})$,  $(\ref{tvd})_1$ is a strongly elliptic system. If  $F \in W^{-1,2}(\Omega)$, then there exists a unique weak solution  $U \in D^1_0(\Omega)$. We begin with recalling various estimates for this system in $L^l(\Omega)$ spaces, which can be seen in \cite{dd8}.

\begin{lemma}\label{tvd1} Let $l \in (1, +\infty)$ and $u$ be a solution of (\ref{tvd}). There exists a constant $C$ depending only on $\lambda$, $\mu$, $l$ and $\Omega$ such that the following estimates hold:

(1) if $F\in L^l(\Omega)$, then we have
\begin{equation}\label{tvd2}
\|U\|_{W^{2,l}}\leq C|F|_l;
\end{equation}

(2)
if $F \in W^{-1,l}(\Omega)$ (i.e., $F =\text{div}f$ with $f =(f_{ij})_{3\times3}$,  $f_{ij} \in L^l(\Omega)$), then we have
\begin{equation}\label{tvd3}
\| U \|_{W^{1,l}} \leq C|f|_l;
\end{equation}

 (3) if $F = \text{div}f$ with $f_{ij} = \partial_k h^k_{ij}$ and $h^k_{ij} \in W^{1,l}_0(\Omega)$ for $ i,j,k = 1,2,3$, then we have
\begin{equation}\label{tvd4}
| U |_{l} \leq C|h|_l.
\end{equation}
\end{lemma}

Moreover, we  need an endpoint estimate for $L$ in the case $l = \infty$. Let $BMO(\Omega)$ stands for the John-Nirenberg space of
bounded mean oscillation whose norm is defined by:
\begin{equation}\label{tvd5}
\|F\|_{BMO(\Omega)}=\|f\|_{L^2(\Omega)}+[f]_{[BMO]},
\end{equation}
with
\begin{equation}\label{tvd6}
\begin{cases}
\displaystyle
[f]_{[BMO]}=\sup_{x\in \Omega,\ r\in (0,d)} \frac{1}{|\Omega_r(x)|}\int_{\Omega_r(x)} |f(y)-f_{\Omega_r(x)}|\text{d}y,\\[10pt]
\displaystyle
f_{\Omega_r(x)}= \frac{1}{|\Omega_r(x)|} \int_{\Omega_r(x)} f(y)\text{d}y,
\end{cases}
\end{equation}
where  $\Omega_{r} (x) = B_r(x)\cap \Omega$, $B_r(x)$ is the ball with center $x$ and radius $r$, and $d$ is the diameter of $\Omega$. $|\Omega_r(x)|$ denotes the Lebesgue measure of $\Omega_r(x)$. Note that
\begin{equation}\label{tvd77}
\begin{split}
[f]_{[BMO]}\leq 2|f|_\infty.
\end{split}
\end{equation}

\begin{lemma}\label{tvd18}
If $F = \text{div} f$ with $f = (f_{ij} )_{3\times3}, \ f_{ij} \in  L^\infty(\Omega)\cap L^2(\Omega)$, then $\nabla U \in  BMO(\Omega)$ and there exists a
constant $C$ depending only on $\lambda$, $\mu$ and $\Omega$ such that
\begin{equation}\label{tvd8}
\begin{split}
|\nabla U|_{[BMO]}\leq C(|f|_\infty+|f|_2).
\end{split}
\end{equation}
\end{lemma}
Because $\Omega$ is a bounded domain with smooth boundary, the estimate (\ref{tvd8}) can be found in \cite{dd8} for a more general setting.
In the next lemma, we will give a variant of the Brezis-Waigner inequality \cite{dd7}, which also can be seen in \cite{zif}.
\begin{lemma}\cite{dd7}\label{tvd10} Let $\Omega$ be a bounded Lipschitz domain and $f \in W^{ 1,l} (\Omega)$ with $l \in  (3, \infty)$. There exists a constant $C$ depending on $l$ and the Lipschitz property of $\Omega$ such that
\begin{equation}\label{tvd9}
\begin{split}
|f|_{L^\infty(\Omega)}\leq& C\big(1+|f|_{BMO(\Omega)} \ln (e+|\nabla  f|_l\big).
\end{split}
\end{equation}
\end{lemma}

Finally, for $(H,u) \in C^1(\Omega)$,  there are some formulas based on  $\text{div}H=0$:
\begin{lemma}\label{liu6} Let $(H, \rho, u,P)$ be the unique strong solution obtained in Theorem \ref{th5} to   IBVP (\ref{eq:1.2})--(\ref{fan1}) in $[0,T)\times \Omega$, then we have
\begin{equation}\label{zhoumou}
\begin{cases}
\displaystyle
\text{rot}(u\times H)=(H\cdot \nabla)u-(u\cdot \nabla)H-H\text{div}u,\\[8pt]
\displaystyle
\text{rot}H\times H
=\text{div}\Big(H\otimes H-\frac{1}{2}|H|^2\mathbb{I}_3\Big)=-\frac{1}{2}\nabla |H|^2+H\cdot \nabla H.
\end{cases}
\end{equation}
\end{lemma}
\begin{proof} It follows immediately from the following equality:
\begin{equation*}
\begin{cases}
a\times \text{rot}a=\frac{1}{2}\nabla (| a|^2)-a \cdot \nabla a,\\[8pt]
\text{rot}(a\times b)=(b\cdot \nabla)a-(a \cdot \nabla)b+(\text{div}b)a-(\text{div}a)b,
\end{cases}
\end{equation*}
based on the fact that $\text{div}H=0$.
\end{proof}

\section{Blow-up criterion (\ref{eq:2.91}) for strong solutions}

 Now we prove (\ref{eq:2.91}). Let $(H, \rho, u,P)$ be the unique strong solution obtained in Theorem \ref{th5} to   IBVP (\ref{eq:1.2})--(\ref{fan1}) in $[0,\overline{T})\times \Omega$.
Due to $P=A\rho^\gamma$, we show that $P$  satisfies
\begin{equation}\label{mou9}
P_t+u\cdot \nabla P+\gamma P \text{div}u=0, \quad P_0 \in H^2 \cap W^{2,q}.
\end{equation}

We first give  the standard energy estimate that

  \begin{lemma}\label{s2}
\begin{equation*}
\begin{split}
|\sqrt{\rho}u(t)|^2_{ 2}+|H(t)|^2_2+|P(t)|_1+\int_{0}^{T}|\nabla u(t)|^2_{2}\text{d}t\leq C,\quad \text{for} \quad 0\leq t\leq  T,
\end{split}
\end{equation*}
where $C$ only depends on $C_0$, $\mu$, $\lambda$, $A$, $\gamma$, $\Omega$ and $T$ $(any\  T\in (0,\overline{T}])$.
 \end{lemma}
\begin{proof} We first show that
\begin{align}
\label{2}\frac{d}{dt}\int_{\Omega} \Big(\frac{1}{2}\rho |u|^2+\frac{P}{\gamma-1}+\frac{1}{2}H^2\Big) \text{d}x+\int_{\Omega} \big(\mu |\nabla u|^2+(\lambda+\mu)(\text{div}u)^2\big)\text{d}x=0.
\end{align}

Actually,  (\ref{2})  is classical, which can be shown by multiplying   $(\ref{eq:1.2})_4$ by $u$,  $(\ref{eq:1.2})_3$ by $\frac{|u|^2}{2}$  and  $(\ref{eq:1.2})_1$ by $H$, then summing them together and integrating the resulting equation over $\Omega$ by parts, where we have used the fact
\begin{equation}\label{zhu1}
\begin{split}
\int_{\Omega} \text{rot}H \times H \cdot u \text{d}x=\int_{\Omega} -\text{rot}(u \times H)  \cdot H        \text{d}x.
\end{split}
\end{equation}
\end{proof}

Our proof for Theorem \ref{th3} is carried out by contradiction.
We assume that the opposite of (\ref{eq:2.91}) holds, i.e.,
\begin{equation}\label{we11*}
\begin{split}
\lim \sup_{T\mapsto \overline{T}}\Big( |\rho|_{L^\infty([0,T]\times \Omega)}+ |H|_{L^\infty([0,T]\times \Omega)}\Big)=C_0<\infty.
\end{split}
\end{equation}
Based on (\ref{we11*}), we state:
\begin{theorem} \label{extension}
Let  the assumptions shown in Theorem \ref{th3} hold. If we have (\ref{we11*}),
then the strong  solution $(H,\rho,u,P)$  can be continued after the time $\overline{T}$.
\end{theorem}
The rest of this section is devoted to its proof.

First of all, via (\ref{we11*}), we can  improve the energy estimate obtained in Lemma \ref{s2}.
  \begin{lemma}\label{abs3}
If (\ref{ass}) holds,
 then there exists $r\in \big(3,\frac{7}{2}\big)$ such that
\begin{equation}\label{keyq}
\begin{split}
\int_{\Omega}\rho|u(t)|^r \text{dx}\leq C,\quad \text{for} \quad 0\leq t\leq  T,
\end{split}
\end{equation}
 where $C$ only depends on $C_0$, $\mu$, $\lambda$, $A$, $\gamma$, $\Omega$ and $T$ $(any\  T\in (0,\overline{T}])$.
 \end{lemma}
\begin{proof}
 For any $\lambda$ satisfying  (\ref{ass}), there must exists a sufficiently small  constant $\alpha_\lambda>0$:
\begin{equation}\label{qudai}
 3\lambda <(29-\alpha_\lambda) \mu.
\end{equation}
 So we only need to show that (\ref{keyq}) holds under the assumption (\ref{qudai}).

 First, multiplying $ (\ref{eq:1.2})_4$ by $r|u|^{r-2}u$ $(r\geq 3)$ and integrating  over $\Omega$, we have
\begin{equation}\label{lz1}
\begin{split}
& \frac{d}{dt}\int_{\Omega} \rho |u|^r\text{d}x+\int_{\Omega}H_r\text{d}x\\
=&-r(r-2)(\mu+\lambda)\int_{\Omega} \text{div}u |u|^{r-3}u\cdot \nabla |u|\text{d}x\\
&+\int_{\Omega} r P\text{div }(|u|^{r-2}u)\text{d}x-\int_{\Omega} r\big(H\otimes H-\frac{1}{2}|H|^2I_3\big): \nabla (|u|^{r-2}u)\text{d}x,
\end{split}
\end{equation}
where
$$
H_r=r|u|^{r-2}\big(\mu|\nabla u|^2+(\mu+\lambda)|\text{div}u|^2+\mu(r-2)\big| \nabla  |u| \big|^2\big).
$$
For any given $\epsilon_1\in (0,1)$ and $\epsilon_0\in (0,\frac{1}{4})$, we define a nonnegative function which will be determined in $\textbf{Step}\ 2$ as follows
$$
\phi(\epsilon_0,\epsilon_1,r)=\left\{ \begin{array}{llll}
\frac{\mu \epsilon_1(r-1)}{3\big(-\frac{(4-\epsilon_0)\mu}{3}-\lambda+\frac{r^2(\lambda+\mu)}{4(r-1)}\big)}, \quad \text{if}\quad l(\mu,\lambda,r,\epsilon_0)=\frac{r^2(\mu+\lambda)}{4(r-1)} -\frac{(4-\epsilon_0)\mu}{3}-\lambda>0, \\[12pt]
\displaystyle
0,\quad \text{otherwise}. \end{array}\right.
$$

$\textbf{Step}\ 1$: we assume that
\begin{equation}\label{ghu}
\begin{split}
&\int_{\Omega \cap |u|>0}  |u|^r \Big| \nabla \Big(\frac{u}{|u|}\Big)\Big|^2\text{d}x> \phi(\epsilon_0,\epsilon_1,r)\int_{\Omega \cap |u|>0}  |u|^{r-2} \big| \nabla  |u|\big|^2\text{d}x.
\end{split}
\end{equation}
 A direct calculation gives for $|u|>0$:
\begin{equation}\label{popo}
\begin{split}
|\nabla u|^2=|u|^2\Big| \nabla \Big(\frac{u}{|u|}\Big)\Big|^2+\big| \nabla  |u|\big|^2,
\end{split}
\end{equation}
which plays an important role in the proof.

Then by (\ref{lz1}) and  the Cauchy's inequality, we have
\begin{equation}\label{lz3}
\begin{split}
& \frac{d}{dt}\int_{\Omega} \rho |u|^r\text{d}x+\int_{\Omega \cap \{|u|>0\}}H_r\text{d}x\\
=&-r(r-2)(\mu+\lambda)\int_{\Omega \cap \{|u|>0\}} \text{div}u |u|^{\frac{r-2}{2}} |u|^{\frac{r-4}{2}}u\cdot \nabla |u|\text{d}x\\
&+\int_{\Omega} r P\text{div }(|u|^{r-2}u)\text{d}x-\int_{\Omega} r\big(H\otimes H-\frac{1}{2}|H|^2I_3\big): \nabla (|u|^{r-2}u)\text{d}x\\
\leq & r(\mu+\lambda)\int_{\Omega \cap \{|u|>0\}}  |u|^{r-2} |\text{div}u|^2\text{d}x\\
&+\frac{r(r-2)^2(\mu+\lambda)}{4}\int_{\Omega \cap \{|u|>0\}}  |u|^{r-2} \big| \nabla  |u| \big|^2\text{d}x\\
&+\int_{\Omega} r P\text{div }(|u|^{r-2}u)\text{d}x-\int_{\Omega} r\big(H\otimes H-\frac{1}{2}|H|^2I_3\big): \nabla (|u|^{r-2}u)\text{d}x.
\end{split}
\end{equation}
Via H\"older's inequality, Gagliardo-Nirenberg inequality and Young's inequality, we have
\begin{equation}\label{zhu2s}
\begin{split}
J_1=&\int_{\Omega} r P\text{div }(|u|^{r-2}u)\text{d}x\\
\leq&Cr(r-1) \Big(\int_{\Omega}|u|^{r-2}|\nabla u|^2\text{d}x\Big)^{\frac{1}{2}}\Big(\int_{\Omega}|u|^{r-2}P^2\text{d}x\Big)^{\frac{1}{2}}\\
\leq & Cr(r-1)|P|_{\frac{12r}{4r+4}}\Big( \int_{\Omega}|u|^{r-2}|\nabla u|^2\text{d}x\Big)^{\frac{1}{2}}\Big(\int_{\Omega}\big(|u|^{\frac{r}{2}}\big)^6\text{d}x\Big)^{\frac{2(r-2)}{12r}}\\
\leq & Cr(r-1)\Big( \int_{\Omega}|u|^{r-2}|\nabla u|^2\text{d}x\Big)^{\frac{1}{2}}\Big(\int_{\Omega}\big( \nabla |u|^{\frac{r}{2}}\big)^2\text{d}x\Big)^{\frac{(r-2)}{2r}}\\
\leq& \frac{1}{2}\mu r \epsilon_0 \int_{\Omega}|u|^{r-2}|\nabla u|^2\text{d}x+C(\mu, r,\epsilon_0),\\
J_2=&-\int_{\Omega} r\big(H\otimes H-\frac{1}{2}|H|^2I_3\big): \nabla (|u|^{r-2}u)\text{d}x\\
\leq& Cr(r-1)\Big(\int_{\Omega}|u|^{r-2}|\nabla u|^2\text{d}x\Big)^{\frac{1}{2}}\Big(\int_{\Omega}|u|^{r-2}|H|^4\text{d}x\Big)^{\frac{1}{2}}\\
\leq &Cr(r-1)|H^2|_{\frac{12r}{4r+4}}\Big( \int_{\Omega}|u|^{r-2}|\nabla u|^2\text{d}x\Big)^{\frac{1}{2}}\Big(\int_{\Omega}\big(|u|^{\frac{r}{2}}\big)^6\text{d}x\Big)^{\frac{2(r-2)}{12r}}\\
\leq & Cr(r-1)\Big( \int_{\Omega}|u|^{r-2}|\nabla u|^2\text{d}x\Big)^{\frac{1}{2}}\Big(\int_{\Omega}\big( \nabla |u|^{\frac{r}{2}}\big)^2\text{d}x\Big)^{\frac{(r-2)}{2r}}\\
\leq& \frac{1}{2} \mu r \epsilon_0 \int_{\Omega}|u|^{r-2}|\nabla u|^2\text{d}x+C(\mu, r,\epsilon_0),
\end{split}
\end{equation}
where $\epsilon_0\in (0,\frac{1}{4})$  is independent of $r$.
Then combining (\ref{popo})-(\ref{zhu2s}), we quickly have
\begin{equation}\label{lz3kk}
\begin{split}
& \frac{d}{dt}\int_{\Omega} \rho |u|^r\text{d}x+\int_{\Omega \cap \{|u|>0\}}  \mu r(1-\epsilon_0)|u|^{r-2}|\nabla |u| |^2\text{d}x\\
&+\int_{\Omega \cap \{|u|>0\}}  \mu r(1-\epsilon_0)|u|^{r}\Big| \nabla \Big(\frac{u}{|u|}\Big)\Big|^2\text{d}x+\int_{\Omega \cap \{|u|>0\}} \mu r(r-2)\big| \nabla  |u| \big|^2\text{d}x\\
\leq &
\frac{r(r-2)^2(\mu+\lambda)}{4}\int_{\Omega \cap \{|u|>0\}}  |u|^{r-2} \big| \nabla  |u| \big|^2\text{d}x+C(\mu, r,\epsilon_0).
\end{split}
\end{equation}
So according to (\ref{ghu}) and (\ref{lz3kk}), we obtain that
\begin{equation}\label{lz4}
\begin{split}
& \frac{d}{dt}\int_{\Omega} \rho |u|^r\text{d}x+rf(\epsilon_0,\epsilon_1, r)\int_{\Omega \cap \{|u|>0\}}|u|^{r-2}|\nabla |u||^2\text{d}x
\leq C(\mu, r,\epsilon_0),
\end{split}
\end{equation}
where
\begin{equation}\label{xuchen}
\begin{split}
f(\epsilon_0,\epsilon_1, r)=\mu (1-\epsilon_0)\phi(\epsilon_0,\epsilon_1,r)+\mu(r-1-\epsilon_0)-\frac{(r-2)^2(\mu+\lambda)}{4}.
\end{split}
\end{equation}

For  $\phi(\epsilon_0,\epsilon_1,r)$, we need to consider the following two cases: one is  $ l(\mu,\lambda,r,\epsilon_0)>0$ for all $[3,+\infty)$, and the other is $ l(\mu,\lambda,r,\epsilon_0)=0$ for some $r\in [3,+\infty)$.  Because $l(\mu,\lambda,r,\epsilon_0)$ is a monotonic increasing function with respect to $r$, we only need to deal with the cases  $l(\mu,\lambda,3,\epsilon_0)>0$ or $l(\mu,\lambda,3,\epsilon_0)\leq 0$. 

\textbf{Subcase} 1:  $l(\mu,\lambda,3,\epsilon_0)>0$, which means that $(5-8\epsilon_0)\mu<3\lambda$ and $ l(\mu,\lambda,r,\epsilon_0)>0$ for all $[3,+\infty)$.
 Therefore, we have
\begin{equation}\label{lzyue}
\begin{split}
\phi(\epsilon_0,\epsilon_1,r)=
\frac{\mu \epsilon_1(r-1)}{3\big(-\frac{(4-\epsilon_0)\mu}{3}-\lambda+\frac{r^2(\lambda+\mu)}{4(r-1)}\big)}
\end{split}
\end{equation}
for any $r\in [3,\infty)$.
Substituting (\ref{lzyue}) into (\ref{xuchen}), for $r\in [3,\infty)$, we have
\begin{equation}\label{xuchendd}
\begin{split}
f(\epsilon_0,\epsilon_1, r)=\frac{\mu^2 \epsilon_1(1-\epsilon_0)(r-1)}{3\big(-\frac{(4-\epsilon_0)\mu}{3}-\lambda+\frac{r^2(\lambda+\mu)}{4(r-1)}\big)}+\mu(r-1-\epsilon_0)-\frac{(r-2)^2(\mu+\lambda)}{4}.
\end{split}
\end{equation}
For $(\epsilon_1,r)=(1,3)$, we have
\begin{equation}\label{dulan}
\begin{split}
f(\epsilon_0, 1,3)=&\frac{16\mu^2(1-\epsilon_0)}{3\lambda-(5-8\epsilon_0)\mu}+\mu(2-\epsilon_0)-\frac{\mu+\lambda}{4}
=-C_1(\lambda-a_1 \mu)(\lambda-a_2\mu),
\end{split}
\end{equation}
then according to $\frac{(5-8\epsilon_0)\mu}{3}<\lambda$, we have $C_1=\frac{3}{4\big(3\lambda-(5-8\epsilon_0)\mu\big)}>0$ and
\begin{equation}\label{dulann}
\begin{split}
a_1(\epsilon_0)=&\frac{13-10\epsilon_0+2\sqrt{64+\epsilon^2_0-56\epsilon_0}}{3},\\
a_2(\epsilon_0)=&\frac{13-10\epsilon_0-2\sqrt{64+\epsilon^2_0-56\epsilon_0}}{3}<0.
\end{split}
\end{equation}
Then if we want to make sure that   $f(\epsilon_0, 1,3)>0$, we have to assume that
\begin{equation}\label{dulann1}
\begin{split}
\frac{(5-8\epsilon_0)\mu}{3}<\lambda< a_1\mu.
\end{split}
\end{equation}

Due to $a_1(0)=\frac{29}{3}$ and $a'_1(\epsilon_0)<0$ for $\epsilon_0\in (0,1/4)$, such a  $\epsilon_0\in (0,1/4)$ exists up to decreasing the value of $\alpha_\lambda$, which can be always done in (\ref{qudai}).

Since $f(\epsilon_0,\epsilon_1, r)$ is continuous w.r.t.  $(\epsilon_1, r)$ over $[0,1]\times [3,+\infty)$, there exists $
\epsilon_1\in  (0,1)$ and $r\in \big(3,\frac{7}{2}\big)$, such that $
f(\epsilon_0,\epsilon_1, r)\geq 0
$,
which, together with (\ref{lz4})-(\ref{xuchen}), implies that
\begin{equation}\label{lz411}
\begin{split}
& \frac{d}{dt}\int_{\Omega} \rho |u|^r\text{d}x \leq C,\quad \text{for \quad some}\quad r\in \big(3,7/2 \big).
\end{split}
\end{equation}

\textbf{Subcase} 2: $l(\mu,\lambda,3,\epsilon_0)\leq 0$, which means that$(5-8\epsilon_0)\mu\geq 3\lambda$. In this case, for $r\in (3,\frac{7}{2})$,  it is easy to get
\begin{equation}\label{peng}
\begin{split}
&r\Big[\mu (1-\epsilon_0)\phi(\epsilon_0,\epsilon_1,r)+\mu(r-1-\epsilon_0)-\frac{(r-2)^2(\mu+\lambda)}{4}\Big]\\
>&3\Big(\frac{7}{4}\mu-\frac{9(\mu+\lambda)}{16}\Big)=3\Big(\frac{19\mu}{16}-\frac{9\lambda}{16}\Big)
\geq 3\Big(\frac{19\mu}{16}-\frac{3(5-8\epsilon_0)\mu}{16}\Big)>\frac{1}{4}\mu,
\end{split}
\end{equation}
which, together with (\ref{lz4})-(\ref{xuchen}),  implies that
\begin{equation}\label{lz422}
\begin{split}
& \frac{d}{dt}\int_{\Omega} \rho |u|^r\text{d}x+\frac{1}{4}\mu\int_{\Omega \cap \{|u|>0\}}|u|^{r-2}\big| \nabla  |u| \big|^2\text{d}x\leq C,\quad \text{for}\quad r\in \big(3, 7/2 \big).
\end{split}
\end{equation}

$\textbf{Step}$ 2 : we assume that
\begin{equation}\label{ghu11}
\begin{split}
&\int_{\Omega \cap |u|>0}  |u|^r \Big| \nabla \Big(\frac{u}{|u|}\Big)\Big|^2\text{d}x\leq  \phi(\epsilon_0,\epsilon_1,r)\int_{\Omega \cap |u|>0}  |u|^{r-2} \big| \nabla  |u|\big|^2\text{d}x.
\end{split}
\end{equation}
A direct calculation gives for $|u|>0$,
\begin{equation}\label{ghu22}
\begin{split}
\text{div}u=|u|\text{div}\Big(\frac{u}{|u|}\Big)+\frac{u\cdot \nabla |u|}{|u|}.
\end{split}
\end{equation}
Then combining (\ref{ghu22}) and (\ref{lz3})-(\ref{zhu2s}), we quickly have
\begin{equation}\label{lz77}
\begin{split}
& \frac{d}{dt}\int_{\Omega} \rho |u|^r\text{d}x+\int_{\Omega \cap \{|u|>0\}}\mu r(1-\epsilon_0)|u|^{r-2}|\nabla u|^2\text{d}x\\
&+\int_{\Omega \cap \{|u|>0\}}r(\lambda+\mu) |u|^{r-2}|\text{div}u|^2\text{d}x+\int_{\Omega \cap \{|u|>0\}}\mu r(r-2)|u|^{r-2}\big| \nabla  |u| \big|^2\text{d}x\\
=&-r(r-2)(\mu+\lambda)\int_{\Omega \cap \{|u|>0\}} \Big( |u|^{r-2} u\cdot \nabla |u| \text{div}\Big(\frac{u}{|u|}\Big)+ |u|^{r-4} |u\cdot \nabla |u| |^2\Big)\text{d}x.
\end{split}
\end{equation}
This gives
\begin{equation}\label{lz88}
\begin{split}
& \frac{d}{dt}\int_{\Omega} \rho |u|^r\text{d}x+\int_{\Omega \cap \{|u|>0\}}r|u|^{r-4}G\text{d}x\leq C(\mu,r,\epsilon_0),
\end{split}
\end{equation}
where
\begin{equation}\label{wang1}
\begin{split}
G=&\mu (1-\epsilon_0)|u|^{2} |\nabla u|^2+(\mu+\lambda)|u|^{2}|\text{div}u|^2+\mu(r-2)|u|^{2}\big| \nabla  |u| \big|^2\\
&+(r-2)(\mu+\lambda)|u|^{2} u\cdot \nabla |u|\text{div}\Big(\frac{u}{|u|}\Big)+(r-2)(\mu+\lambda)|u \cdot \nabla |u||^2.
\end{split}
\end{equation}
Now we consider how to make sure that $G\geq 0$.
\begin{equation}\label{wang2mm}
\begin{split}
G=&\mu (1-\epsilon_0)|u|^{2} \Big( |u|^2\Big| \nabla \Big(\frac{u}{|u|}\Big)\Big|^2+\big| \nabla  |u|\big|^2\Big)+(\mu+\lambda)|u|^{2}|\Big(|u|\text{div}\Big(\frac{u}{|u|}\Big)+\frac{u\cdot \nabla |u|}{|u|}\Big)^2\\
&+\mu(r-2)|u|^{2}\big| \nabla  |u| \big|^2+(r-2)(\mu+\lambda)|u|^{2} u\cdot \nabla |u|\text{div}\Big(\frac{u}{|u|}\Big)\\
&+(r-2)(\mu+\lambda)|u \cdot \nabla |u||^2\\
=&\mu(1-\epsilon_0) |u|^{4} \Big| \nabla \Big(\frac{u}{|u|}\Big)\Big|^2+\mu(r-1-\epsilon_0)|u|^2\big| \nabla  |u| \big|^2+(r-1)(\mu+\lambda)|u \cdot \nabla |u||^2\\
&+r(\mu+\lambda)|u|^{2} u\cdot \nabla |u|\text{div}\Big(\frac{u}{|u|}\Big)+(\mu+\lambda)|u|^{4}\Big(\text{div}\Big(\frac{u}{|u|}\Big)\Big)^2\\
=&\mu(1-\epsilon_0) |u|^{4} \Big| \nabla \Big(\frac{u}{|u|}\Big)\Big|^2+\mu(r-1-\epsilon_0)|u|^{2}\big| \nabla  |u| \big|^2\\
&+(r-1)(\mu+\lambda)\Big(u \cdot \nabla |u|+\frac{r}{2(r-1)}|u|^{2}\Big(\text{div}\frac{u}{|u|}\Big)\Big)^2\\
&+(\mu+\lambda)|u|^{4}\Big(\text{div}\frac{u}{|u|}\Big)^2-\frac{r^2(\mu+\lambda)}{4(r-1)}|u|^{4}\Big(\text{div}\Big(\frac{u}{|u|}\Big)\Big)^2,
\end{split}
\end{equation}
which, combining with the fact
$$
\Big|\text{div}\Big(\frac{u}{|u|}\Big)\Big|^2\leq 3\Big|\nabla \Big(\frac{u}{|u|}\Big)\Big|^2,
$$
implies that
\begin{equation}\label{wang3}
\begin{split}
G\geq& \mu(1-\epsilon_0) |u|^{4} \Big| \nabla \Big(\frac{u}{|u|}\Big)\Big|^2+\mu(r-1-\epsilon_0)|u|^{2}\big| \nabla  |u| \big|^2\\
&+\Big(\mu+\lambda-\frac{r^2(\mu+\lambda)}{4(r-1)}\Big)|u|^{4}\Big(\text{div}\Big(\frac{u}{|u|}\Big)\Big)^2\\
\geq &\frac{\mu(1-\epsilon_0)}{3} |u|^{4} \Big| \text{div}\Big(\frac{u}{|u|}\Big)\Big|^2+\mu(r-1-\epsilon_0)|u|^{2}\big| \nabla  |u| \big|^2\\
&+\Big(\mu+\lambda-\frac{r^2(\mu+\lambda)}{4(r-1)}\Big)|u|^{4}\Big(\text{div}\Big(\frac{u}{|u|}\Big)\Big)^2\\
\geq &\mu(r-1-\epsilon_0)|u|^{2}\big| \nabla  |u| \big|^2+\Big(\frac{(4-\epsilon_0)\mu}{3}+\lambda-\frac{r^2(\mu+\lambda)}{4(r-1)}\Big)|u|^{4}\Big(\text{div}\Big(\frac{u}{|u|}\Big)\Big)^2.
\end{split}
\end{equation}
Thus
\begin{equation}\label{wang4}
\begin{split}
&\int_{\Omega \cap \{|u|>0\}} r|u|^{r-4} G\text{d}x\\
\geq & r \Big(\frac{(4-\epsilon_0)\mu}{3}+\lambda-\frac{r^2(\mu+\lambda)}{4(r-1)}\Big) \int_{\Omega \cap \{|u|>0\}} |u|^{r} \Big(\text{div}\Big(\frac{u}{|u|}\Big)\Big)^2\text{d}x\\
&+\mu r(r-1-\epsilon_0)\int_{\Omega \cap \{|u|>0\}} |u|^{r-2} \big| \nabla  |u| \big|^2\text{d}x\\
\geq & 3r \Big(\frac{(4-\epsilon_0)\mu}{3}+\lambda-\frac{r^2(\mu+\lambda)}{4(r-1)}\Big)\phi(\epsilon_0,\epsilon_1,r)\int_{\Omega \cap \{|u|>0\}}|u|^{r-2} \big| \nabla  |u| \big|^2\text{d}x\\
&+\mu r(r-1-\epsilon_0)\int_{\Omega \cap \{|u|>0\}} |u|^{r-2} \big| \nabla  |u| \big|^2\text{d}x\\
\geq & g(\epsilon_0,\epsilon_1,r)\int_{\Omega \cap \{|u|>0\}}  |u|^{r-2} \big| \nabla  |u| \big|^2\text{d}x,
\end{split}
\end{equation}
where we have used (\ref{ghu11}) and 
\begin{equation}\label{xuchendd}
\begin{split}
g(\epsilon_0,\epsilon_1,r)=\Big[3r \Big(\frac{(4-\epsilon_0)\mu}{3}+\lambda-\frac{r^2(\mu+\lambda)}{4(r-1)}\Big)\phi(\epsilon_0,\epsilon_1,r)+\mu r(r-1-\epsilon_0)\Big)\Big].
\end{split}
\end{equation}
Here we need that $\epsilon_0$ is sufficiently small such that $\epsilon_0<(r-1)(1-\epsilon_1)$.
Then combining (\ref{lz88}) and  (\ref{wang4})-(\ref{xuchendd}), we quickly have
\begin{equation}\label{lz4mm}
\begin{split}
& \frac{d}{dt}\int_{\Omega} \rho |u|^r\text{d}x \leq C,\quad \text{for}\quad r\in \big(3,7/2\big).
\end{split}
\end{equation}

So combining (\ref{lz411})-(\ref{lz422}) and (\ref{lz4mm}) for $\textbf{Step}$:1 and $\textbf{Step}$:2, we conclude that if $3\lambda <(29-\alpha_\lambda)\mu$, there exits some constants $C>0$ such that
\begin{equation}\label{lz4nn}
\begin{split}
& \frac{d}{dt}\int_{\Omega} \rho |u|^r\text{d}x \leq C,\quad \text{for\quad some }\quad r\in \big(3,7/2\big).
\end{split}
\end{equation}

\end{proof}

Now for each $t\in [0,T)$, we denote $v(t,x)=(-L)^{-1} \text{div} A$ and
\begin{equation}\label{dingyiA}
A=PI_3-\Big(H\otimes H-\frac{1}{2}|H|^2I_3\Big),
\end{equation}
that is, $v$ is the solution of
\begin{equation}\label{diyi}
\begin{cases}
\mu \triangle v+(\lambda+\mu)\nabla \text{div} v=\text{div} A \quad \text{in}\quad \Omega,\\[8pt]
v(t,x)=0 \quad \text{on} \quad \partial \Omega.
\end{cases}
\end{equation}
From Lemma \ref{tvd1}, for any $l\in (1,+\infty)$, there exists a constant $C$ independent of $t$ such that
\begin{equation}\label{diy2}
\begin{cases}
|\nabla v(t)|_l\leq C(|\rho(t)|_l+|H(t)|_l),\\[8pt]
|\nabla^2 v(t)|_l\leq C(|\nabla \rho(t)|_l+|\nabla H (t)|_l).
\end{cases}
\end{equation}
Now let us introduce an important quantity:
$$ w=u-v,$$
whose divergence can be viewed as the effective viscous flux (see Hoff \cite{hoff1}).
It will be  observed that this quantity $w$ possesses more regularity information than $u$  under the assumption that $(H,\rho)$ is upper bounded. First, we have
  \begin{lemma}\label{sk4}
\begin{equation*}
\begin{split}
|\nabla w(t)|^2_{ 2}+\int_0^T( |\nabla^2 w|^2_2+|\sqrt{\rho}w_t|^2_2)\text{d}t\leq C,\quad \text{for} \quad 0\leq t\leq  T,
\end{split}
\end{equation*}
 where $C$ only depends on $C_0$, $\mu$, $\lambda$, $A$, $\gamma$, $\Omega$ and $T$ $(any\  T\in (0,\overline{T}])$. \end{lemma}

\begin{proof}

First, from the momentum equations $(\ref{eq:1.2})_4$, we find that $w$ satisfies
\begin{equation}\label{zhu6}
\begin{cases}
\rho w_t -\mu \triangle w-(\lambda+\mu)\nabla \text{div}w=\rho F,\\
w(t,x)=0 \quad \text{on} \quad [0,T)\times \partial \Omega, \ w(0,x)=w_0(x), \quad \text{in }\  \Omega,
\end{cases}
\end{equation}
with $w_0(x)=u_0(x)+v_0(x)$ and
\begin{equation*}
\begin{split}
F=&-u\cdot \nabla u+L^{-1} \text{div} A_t\\
=&-u\cdot \nabla u-L^{-1} \nabla \text{div}(Pu)-(\gamma-1) L^{-1}\nabla (P\text{div}u)\\
&-L^{-1}\text{div}(H_t\otimes H+H\otimes H_t )+L^{-1} \nabla (H\cdot H_t)=\sum_{i=1}^{5} J_i,
\end{split}
\end{equation*}
where the definition of $A$ could be seen in (\ref{dingyiA}).
Multiplying the equations in (\ref{zhu6}) by $w_t$ and integrating the resulting equation over $\Omega$, from H\"older's inequality, we have
\begin{equation}\label{gaibian}
\begin{split}
&\frac{d}{dt}\int_{\Omega} \big(\mu|\nabla w|^2+(\lambda+\mu)|\text{div}w|^2\big)\text{d}x +\int_{\Omega}\rho |w_t|^2 \text{d}x\\
=&\int_{\Omega} \rho F\cdot w_t \text{d}x
 \leq C|\sqrt{\rho}F|^2_2+\frac{1}{2}\int_{\Omega}\rho |w_t|^2 \text{d}x,
\end{split}
\end{equation}
which means that
\begin{equation}\label{gai11}
\begin{split}
&\frac{d}{dt}\int_{\Omega} \big(\mu|\nabla w|^2+(\lambda+\mu)|\text{div}w|^2\big)\text{d}x +\int_{\Omega}\rho |w_t|^2 \text{d}x
\leq  C\sum_{i=1}^{5} |\sqrt{\rho} J_i|^2_2.
\end{split}
\end{equation}
Next we need to consider the terms $|\sqrt{\rho} J_i|_2$ for $i=1,2,...,5$. From Lemma \ref{abs3} and (\ref{diy2}), it follows that
\begin{equation}\label{yue5}
\begin{split}
|\sqrt{\rho}J_1|_2=& |-\sqrt{\rho}u \cdot \nabla u|_2
\leq C |\sqrt{\rho}u |_r |\nabla u|_{\frac{2r}{r-2}}\\
\leq&  C\Big( |\nabla w|_{\frac{2r}{r-2}}+|\nabla v|_{\frac{2r}{r-2}}\Big)
\leq C(\epsilon) |\nabla w|_2+\epsilon | w|_{D^2}+C,
\end{split}
\end{equation}
where we have used the interpolation inequality
$$
|f|_p\leq C(\epsilon)|f|_2+\epsilon |\nabla f|_2,\quad 2\leq p <6.
$$
According to Lemmas \ref{s2}-\ref{abs3}, we obtain
\begin{equation}\label{kkll}
\begin{split}
|\sqrt{\rho}J_2|_2=& |-\sqrt{\rho}L^{-1} \nabla \text{div}(Pu)|_2\leq C|Pu|_2\leq C|\sqrt{\rho}u|_2\leq C,\\
|\sqrt{\rho}J_3|_2=& |-(\gamma-1)\sqrt{\rho}L^{-1} \nabla (P\text{div}u)|_2\\
\leq&C |\sqrt{\rho}|_3|L^{-1} \nabla (P\text{div}u)|_6\\
\leq & C |\nabla L^{-1} \nabla (P\text{div}u)|_2\leq  C |P\text{div}u |_2\leq  C|\nabla u|_2.
\end{split}
\end{equation}
Now we consider the term $\mathbb{B}=(b^{(i,j)})_{(3\times 3)}=H_t\otimes H+H\otimes H_t$.
Due to  Lemma \ref{liu6},
\begin{equation}\label{bang1}
H_t=H \cdot \nabla u-u \cdot \nabla H-H\text{div}u,
\end{equation}
we get
\begin{equation}\label{yue66}
\begin{split}
b^{(i,j)}=&H^j\big(H^k\partial_k u^i-u^k\partial_k H^i-H^i\partial_k u^k\big)\\
&+H^i\big(H^k\partial_k u^j-u^k\partial_k H^j-H^j\partial_k u^k\big)\\
=&H^iH^k\partial_k u^j+H^jH^k\partial_k u^i-H^i H^j\partial_k u^k-\partial_k(H^iH^j u^k),
\end{split}
\end{equation}
which means that
\begin{equation}\label{yue67}
\begin{split}
\mathbb{B}=&(H^iH^k\partial_k u^j+H^jH^k\partial_k u^i-H^i H^j\partial_k u^k)_{(3\times 3)}-\text{div}\big( (H\otimes H) \otimes u\big)
=\mathbb{B}_1+\mathbb{B}_2.
\end{split}
\end{equation}
Then we   have
\begin{equation}\label{kksd}
\begin{split}
|\sqrt{\rho}J_4|_2=& |\sqrt{\rho}L^{-1}\text{div}(H_t\otimes H+H\otimes H_t )|_2\\
=& |\sqrt{\rho}L^{-1} \text{div}\mathbb{B}_1|_2+ |\sqrt{\rho}L^{-1} \text{div}\mathbb{B}_2 |_2\\
\leq&C |\sqrt{\rho}|_3|L^{-1} \text{div}\mathbb{B}_1|_6+|\sqrt{\rho}L^{-1} \text{div}\text{div}\big( (H\otimes H)\otimes u\big) |_2\\
\leq & C |\nabla L^{-1} \text{div}\mathbb{B}_1|_2+C|\nabla u|_2\leq C|\nabla u|_2.
\end{split}
\end{equation}
Similarly, we consider the term $\mathbb{C}=H\cdot H_t$. Due to (\ref{bang1}), we obtain
\begin{equation}\label{yue68}
\begin{split}
\mathbb{C}=& H\cdot (H \cdot \nabla u-u \cdot \nabla H-H\text{div}u)\\
=&\Big(H\cdot \nabla u\cdot H-\frac{1}{2}|H|^2\text{div}u\Big)-\frac{1}{2}\text{div}(u|H|^2)
=\mathbb{C}_1+\mathbb{C}_2,
\end{split}
\end{equation}
which, together with the Poincar$\acute{\text{e}}$ inequality,  implies that
\begin{equation}\label{kks}
\begin{split}
|\sqrt{\rho}J_5|_2=& |\sqrt{\rho}L^{-1} \nabla (H\cdot H_t)|_2\\
=& |\sqrt{\rho}L^{-1} \nabla \mathbb{C}_1|_2+ |\sqrt{\rho}L^{-1} \nabla \mathbb{C}_2 |_2\\
\leq&C |\sqrt{\rho}|_3|L^{-1} \nabla \mathbb{C}_1|_6+|\sqrt{\rho}L^{-1} \nabla \text{div}(u|H|^2) |_2\\
\leq & C |\nabla L^{-1} \nabla \mathbb{C}_1|_2+C|\nabla u|_2\leq C|\nabla u|_2.
\end{split}
\end{equation}
Combining (\ref{yue5})-(\ref{kks}), we have
\begin{equation}\label{kksk}
\begin{split}
|\sqrt{\rho}F|^2_2\leq \epsilon|\nabla^2 w|^2_2+C_1(\epsilon)(1+|\nabla w|^2_2+|\nabla u|^2_2).
\end{split}
\end{equation}
where $C_1(\epsilon)>0$ is a constant dependent of $\epsilon$.
Then from Lemma \ref{tvd1} and (\ref{zhu6}), we have
\begin{equation}\label{kksu}
\begin{split}
|\nabla^2 w|^2_2\leq C(|\rho w_t|^2_2+|\rho F|^2_2)\leq C_2(|\sqrt{\rho} w_t|^2_2+|\sqrt{\rho} F|^2_2),
\end{split}
\end{equation}
where $C_2>0$ is a constant independent of $\epsilon$.
Then we have
\begin{equation}\label{kksw}
|\sqrt{\rho}F|^2_2\leq  C_2\epsilon (|\sqrt{\rho} w_t|^2_2+|\sqrt{\rho} F|^2_2)+C_1(\epsilon)(1+|\nabla w|^2_2+|\nabla u|^2_2).
\end{equation}
Taking $\epsilon=\frac{1}{3C_2}$ in (\ref{kksw}), we obtain
\begin{equation}\label{kksw1}
\begin{split}
|\sqrt{\rho}F|^2_2
\leq& \frac{1}{2}|\sqrt{\rho} w_t|^2_2 +\frac{3}{2}C_1(\epsilon)(1+|\nabla w|^2_2+|\nabla u|^2_2).
\end{split}
\end{equation}

Substituting (\ref{kksw1}) into (\ref{gai11}), from Gronwall's inequality, the desired conclusions can be obtained.
\end{proof}

Finally, according to the estimates obtained in  (\ref{diy2}) and  Lemmas \ref{abs3}-\ref{sk4}, we get
  \begin{lemma}\label{sk4ss}
\begin{equation*}
\begin{split}
|\nabla u(t)|^2_{ 2}+\int_0^T |\nabla u|^2_q\text{d}t\leq C,\quad \text{for} \quad 0\leq t\leq  T, \quad \text{and} \quad q\in (3,6],
\end{split}
\end{equation*}
 where $C$ only depends on $C_0$, $\mu$, $\lambda$, $A$, $\gamma$, $\Omega$ and $T$ $(any\  T\in (0,\overline{T}])$. \end{lemma}

Next, we will give high order regularity estimates for $w$. This is possible if the initial data $(H_0,\rho_0,u_0, P_0)$ satisfies the compatibility condition (\ref{th79}). First for a function or vector field (or even a $3 \times 3$ matrix) $f (t , x )$, the material derivative $ f$  is defined by:
$$
\dot{f}=f_t+u\cdot \nabla f=f_t+\text{div}(fu)-f\text{div}u.
$$

 \begin{lemma}[\textbf{Lower order estimate of the velocity $u$}]\label{ablem:4-1}\ \\
\begin{equation*}
\begin{split}
|w(t)|^2_{D^2}+|\sqrt{\rho} \dot{u}(t)|^2_{2} + \int_{0}^{T} |\dot{u}|^2_{D^1}\text{d}t\leq C,\quad \text{for} \quad 0\leq t\leq  T,
\end{split}
\end{equation*}
 where $C$ only depends on $C_0$, $\mu$, $\lambda$, $A$, $\gamma$, $\Omega$ and $T$ $(any\  T\in (0,\overline{T}])$.
 \end{lemma}

 \begin{proof}
We will follow an idea due to Hoff \cite{hoff1}. Applying $\dot{u}[\partial / \partial t+\text{div}(u\cdot)]$ to $(\ref{eq:1.2})_4$  and integrating by parts give
\begin{equation}\label{bzhen4}
\begin{split}
&\frac{1}{2}\frac{d}{dt}\int_{\Omega}\rho |\dot{u}|^2 \text{d}x\\
=&  -\int_{\Omega}\Big(   \dot{u}\cdot \big(\nabla P_t+\text{div}(\nabla P\otimes u)\big)+\dot{u}\cdot \big(\triangle u_t+\text{div}(\triangle u\otimes u)\big)  \Big)\text{d}x\\
&+ (\lambda+\mu)\int_{\Omega} \dot{u}\cdot  \big(\nabla \text{div}u_t+\text{div}(\nabla  \text{div}u\otimes u)\big) \text{d}x\\
&+\int_{\Omega} \dot{u} \cdot \Big(\text{div}\Big(H\otimes H-\frac{1}{2}|H|^2I_3\Big)_t+\text{div}\Big(\text{div}\Big(H\otimes H-\frac{1}{2}|H|^2I_3\Big)\otimes u\Big)\Big)\text{d}x\\
\equiv&:\sum_{i=6}^{8}J_i+\Lambda.
\end{split}
\end{equation}

According to Lemmas \ref{s2}-\ref{sk4ss}, equation (\ref{mou9}),  H\"older's inequality, Gagliardo-Nirenberg inequality and Young's inequality,  we deduce that
\begin{equation}\label{zhou6}
\begin{split}
J_6=& -\int_{\Omega}\big(\dot{u}\cdot \big(\nabla P_t+\text{div}(\nabla P\otimes u)\big)\big)\text{d}x\\
=& \int_{\Omega}\big(\partial_j\dot{u}^j P_t+\partial_k\dot{u}^j\partial_j Pu^k\big)\text{d}x\\
=& \int_{\Omega}\big(-\partial_j \dot{u}^ju^k\partial_k P-\gamma P\text{div}u\partial_j\dot{u}^j +\partial_k\dot{u}^j\partial_j Pu^k\big)\text{d}x\\
=& \int_{\Omega}\big(-\gamma P\text{div}u\partial_j\dot{u}^j +P\partial_k(\partial_j \dot{u}^ju^k)-P\partial_j(\partial_k \dot{u}^ju^k)\big)\text{d}x\\
\leq&  C|\nabla \dot{u}|_{2}|\nabla u|_{2}
\leq \epsilon |\nabla \dot{u}|^2_{2} +C(\epsilon) |\nabla u|^2_2,\\
J_7=& \int_{\Omega}\mu\big(\dot{u}\cdot \big(\triangle u_t+\text{div}(\triangle u\otimes u)\big)\big)\text{d}x\\
=& -\int_{\Omega}\mu\big(\partial_i\dot{u}^j\partial_i u^j_t+\triangle u^j u\cdot \nabla \dot{u}^j\big)\text{d}x\\
=& -\int_{\Omega}\mu\big(|\nabla \dot{u}|^2- \partial_i \dot{u}^j u^k\partial_k \partial_i u^j- \partial_i \dot{u}^j \partial_i u^k\partial_k u^j+\triangle u^j u\cdot \nabla \dot{u}^j\big)\text{d}x\\
=& -\int_{\Omega}\mu\big(|\nabla \dot{u}|^2+ \partial_i \dot{u}^j \partial_ku^k \partial_i u^j- \partial_i \dot{u}^j \partial_i u^k\partial_k u^j-\partial_iu^j \partial_i u^k \partial_k \dot{u}^j\big)\text{d}x\\
\leq &-\frac{\mu}{2}|\nabla \dot{u}|^2_{2}+C|\nabla u|^4_4,
\end{split}
\end{equation}
where $\epsilon>0$ is a sufficiently small constant.
Similarly, we have
\begin{equation}\label{yueyue}
\begin{split}
J_8=&  (\lambda+\mu)\int_{\Omega} \big(\dot{u}\cdot  \big(\nabla \text{div}u_t+\text{div}(\nabla  \text{div}u\otimes u\big)\big)\text{d}x
\leq  -\frac{\mu+\lambda}{2}|\nabla \dot{u}|^2_2 +C |\nabla u|^4_4.
\end{split}
\end{equation}
Next we begin to consider the magnetic term $\Lambda$
\begin{equation*}
\begin{split}
\Lambda=& \int_{\Omega} \dot{u} \cdot \Big(\text{div}\Big(H\otimes H-\frac{1}{2}|H|^2I_3\Big)_t+\text{div}\Big(\text{div}\Big(H\otimes H-\frac{1}{2}|H|^2I_3\Big)\otimes u\Big)\Big)\text{d}x\
= \sum_{j=1}^4 \Lambda_j.
\end{split}
\end{equation*}
Via the  magnetic equations  $(\ref{eq:1.2})_1$  and integrating by parts, we obtain that
\begin{equation}\label{yueyue12}
\begin{split}
\Lambda_1=& \int_{\Omega} \dot{u} \cdot \text{div}\big(H\otimes H\big)_t\text{d}x= -\int_{\Omega} \big(H\otimes H\big)_t : \nabla \dot{u} \text{d}x\\
=&- \int_{\Omega} \big(H\otimes H_t+H_t\otimes H): \nabla \dot{u} \text{d}x\\
=& -\int_{\Omega}H \otimes\big( H\cdot \nabla u-u \cdot \nabla H-H\text{div}u\big): \nabla \dot{u}\text{d}x\\
&- \int_{\Omega} \big(H \cdot \nabla u-u \cdot \nabla H-H\text{div}u\big)\otimes H: \nabla \dot{u}\text{d}x\\
=& -\int_{\Omega}H \otimes\big( H\cdot \nabla u-H\text{div}u\big): \nabla \dot{u}\text{d}x- \int_{\Omega} \big(H \cdot \nabla u-H\text{div}u\big)\otimes H: \nabla \dot{u}\text{d}x\\
&-\int_{\Omega} \Big(-H\otimes(u \cdot \nabla H)-(u \cdot \nabla H)\otimes H)\Big): \nabla \dot{u}\text{d}x=\Lambda_{11}+\Lambda_{12}+\Lambda_{13},\\
\end{split}
\end{equation}
\begin{equation}\label{yueyue13}
\begin{split}
\Lambda_2=& \int_{\Omega} \dot{u} \cdot \text{div}\Big(-\frac{1}{2}|H|^2I_3\Big)_t\text{d}x= \int_{\Omega} \Big(\frac{1}{2}|H|^2I_3\Big)_t : \nabla \dot{u} \text{d}x\\
=& \int_{\Omega} (H\cdot H_t I_3): \nabla \dot{u} \text{d}x\\
=& \int_{\Omega}\big(H \cdot ( H\cdot \nabla u-u \cdot \nabla H-H\text{div}u)I_3\big): \nabla \dot{u}\text{d}x\\
=& \int_{\Omega}\big(H \cdot\big( H\cdot \nabla u-H\text{div}u\big)I_3\big): \nabla \dot{u}\text{d}x-\int_{\Omega}\big(H \cdot\big( u \cdot \nabla H\big)I_3\big): \nabla \dot{u}\text{d}x\\
=&\Lambda_{21}+\Lambda_{22},\\
\Lambda_3=& \int_{\Omega} \dot{u} \cdot \text{div}\big(\text{div}   \big(H\otimes H\big) \otimes u\big) \text{d}x
= -\int_{\Omega} \text{div}   \big(H\otimes H\big) \otimes u: \nabla  \dot{u}\text{d}x\\
=&- \int_{\Omega}   (H\cdot \nabla H)\otimes u: \nabla  \dot{u}\text{d}x
= -\int_{\Omega}   H^k \partial_kH^i u^j \partial_j \dot{u}^i\text{d}x\\
=& \int_{\Omega}  H^k H^i  \partial_k u^j\partial_j\dot{u}^i\text{d}x + \int_{\Omega}  H^k H^i u^j  \partial_{kj}\dot{u}^i\text{d}x =\Lambda_{31}+\Lambda_{32},\\
\Lambda_4=& \int_{\Omega} \dot{u} \cdot \text{div}\Big(\text{div}\Big(-\frac{1}{2}|H|^2I_3\Big)\otimes u\Big)\text{d}x\\
=&\int_{\Omega}\text{div}\Big(\frac{1}{2}|H|^2I_3\Big)\otimes u:\nabla \dot{u} \text{d}x
=\int_{\Omega}H^k\partial_iH^ku^j\partial _j \dot{u}^i \text{d}x\\
=&-\frac{1}{2}\int_{\Omega} |H^k|^2 \partial_iu^j\partial _j\dot{u}^i \text{d}x-\frac{1}{2}\int_{\Omega} |H^k|^2 u^j\partial _{ij}\dot{u}^i \text{d}x=\Lambda_{41}+\Lambda_{42},
\end{split}
\end{equation}
where we have used the fact that $\text{div}H=0$. Now we observe that
\begin{equation}\label{yueyue15}
\begin{split}
\Lambda_{13}+\Lambda_{32}
=&\int_{\Omega} \Big(H\otimes(u \cdot \nabla H)+(u \cdot \nabla H)\otimes H)\Big): \nabla \dot{u}\text{d}x+ \int_{\Omega}  H^k H^i u^j  \partial_{kj}\dot{u}^i\text{d}x\\
=& \int_{\Omega}  \Big( H^i u^k \partial_{k}H^j  \partial_{j}\dot{u}^i+u^k \partial_{k}H^i H^j \partial_{j}\dot{u}^i+H^k H^i u^j  \partial_{kj}\dot{u}^i\Big)\text{d}x\\
=& \int_{\Omega}  \Big(-\partial_{k}u^k H^i H^j  \partial_{j}\dot{u}^i-u^k H^j\partial_{k}H^i  \partial_{j}\dot{u}^i-H^j H^i u^k \partial_{kj}\dot{u}^i\Big)\text{d}x\\
&+\int_{\Omega}  \Big(u^k H^j\partial_{k}H^i  \partial_{j}\dot{u}^i+H^k H^i u^j  \partial_{kj}\dot{u}^i\Big)\text{d}x\\
=& -\int_{\Omega}  \partial_{k}u^k H^i H^j  \partial_{j}  \dot{u}^i  \text{d}x,\\
\Lambda_{22}+\Lambda_{42}
=&-\int_{\Omega}\big(H \cdot\big( u \cdot \nabla H\big)I_3\big): \nabla \dot{u}\text{d}x-\frac{1}{2}\int_{\Omega} |H^k|^2 u^j\partial _{ij}\dot{u}^i \text{d}x\\
=&-\int_{\Omega} \Big( H^ku^l \partial_l H^k \text{div}\dot{u}+ \frac{1}{2} |H^k|^2 u^j \partial _{ij} \dot{u}^i\Big) \text{d}x\\
=&\int_{\Omega} \Big(  \frac{1}{2}u^j |H^k|^2 \partial_{ij}\dot{u}^i+\frac{1}{2}\partial_ju^j |H^k|^2 \partial_{i}\dot{u}^i-\frac{1}{2} |H^k|^2 u^j \partial _{ij} \dot{u}^i\Big) \text{d}x\\
=&\frac{1}{2}\int_{\Omega} \partial_ju^j |H^k|^2 \partial_{i}\dot{u}^i \text{d}x,
\end{split}
\end{equation}
which, together with (\ref{yueyue12})-(\ref{yueyue13}),  implies that
\begin{equation}\label{yueyue18}
\begin{split}
\Lambda  \leq C|H|^2_{\infty} |\nabla u|_2|\nabla \dot{u}|_2\leq \epsilon|\nabla \dot{u}|^2_2+  C(\epsilon)|\nabla u|^2_2.
\end{split}
\end{equation}
Due to the definition of $w$, we know that $w$ satisfies
\begin{equation}\label{mou5k}
\begin{split}
\mu \triangle w+(\lambda+\mu)\nabla \text{div}w=\rho \dot{u}\quad \text{in} \ \Omega,
\end{split}
\end{equation}
with the zero boundary condition. From Lemma \ref{tvd1}, we have
\begin{equation}\label{mou5s}
|w|_{D^2}\leq C|\rho \dot{u}|_2\leq C|\rho \dot{u}|_2,
\end{equation}
which, together with (\ref{bzhen4})-(\ref{yueyue18}) and  letting $\epsilon>0$ be sufficiently small,  implies that
\begin{equation}\label{mou5}
\begin{split}
\frac{1}{2}\frac{d}{dt}  & \int_{\Omega}\rho |\dot{u}|^2 \text{d}x+|\dot{u}|^2_{D^1}
\leq  C|\nabla u|^4_4+C\\
\leq & C(|\nabla u|_2  |\nabla u|^3_6+1)\leq C\big(|\nabla u|^2_6(|\nabla w|_6+|\nabla v|_6)+1\big)\\
\leq  & C\big(|\nabla u|^2_6(1+|\nabla^2 w|_2)+1\big)
\leq  C\big(|\nabla u|^2_6(1+|\sqrt{\rho}\dot{u}|_2)+1\big).
\end{split}
\end{equation}

Integrating (\ref{mou5}) over $(\tau,t)$ ($\tau\in (0,t)$), for $\tau\leq t \leq T$, we  have
\begin{equation}\label{nvk3}
\begin{split}
&|\sqrt{\rho}\dot{u}(t)|^2_2 +\int_{\tau}^{t}|\nabla \dot{u}|^2_{2}\text{d}s
\leq |\sqrt{\rho}\dot{u}(\tau)|^2_2+C\int_{\tau}^{t} |\nabla u|^2_6(1+|\sqrt{\rho}\dot{u}|_2) \text{d}s+C.
\end{split}
\end{equation}
From the momentum equations  $(\ref{eq:1.2})_3$ , we easily have
\begin{equation}\label{li9}
\begin{split}
|\sqrt{\rho}\dot{u}(\tau)|^2_2\leq C\int_{\Omega} \frac{|\nabla P+Lu- \text{rot}H\times H|^2}{\rho}(\tau)\text{d}x.
\end{split}
\end{equation}
Due to the initial layer compatibility condition (\ref{th79}), letting $\tau\rightarrow 0$ in (\ref{li9}), we have
\begin{equation}\label{nvk33}
\begin{split}
\lim \sup_{\tau\rightarrow 0}|\sqrt{\rho}\dot{u}(\tau)|^2_2 \leq C\int_{\Omega} |g_1|^2\text{d}x\leq C.
\end{split}
\end{equation}
Then, letting $\tau\rightarrow 0$ in (\ref{nvk3}), we have
\begin{equation}\label{nvk44}
\begin{split}
&|\sqrt{\rho}\dot{u}(t)|^2_2 +\int_{0}^{t}|\nabla \dot{u}|^2_{2}\text{d}s \leq C\int_{0}^{t} |\nabla u|^2_6(1+|\sqrt{\rho}\dot{u}|_2) \text{d}s+C.
\end{split}
\end{equation}

Then from Gronwall's inequality and Lemma \ref{sk4ss}, we have
\begin{equation}\label{mou999}
\begin{split}
 \int_{\Omega}\rho |\dot{u}|^2(t) \text{d}x+\int_0^t|\dot{u}|^2_{D^1}\leq & C,\quad \text{for} \quad 0\leq t\leq  T.
\end{split}
\end{equation}
\end{proof}

According to Lemmas \ref{abs3}-\ref{ablem:4-1} and using the equations (\ref{mou5k}) again, we deduce

  \begin{lemma}\label{sk4zz}
\begin{equation*}
\begin{split}
|\nabla w|_{ L^2([0,T];L^{\infty}(\Omega))}+|\nabla^2 w|_{ L^2([0,T];L^{q}(\Omega))}\leq C,\quad \text{for} \quad 0\leq t\leq  T, \quad \text{and} \quad q\in (3,6],
\end{split}
\end{equation*}
 where $C$ only depends on $C_0$, $\mu$, $\lambda$, $A$, $\gamma$, $\Omega$ and $T$ $(any\  T\in (0,\overline{T}])$. \end{lemma}

Finally, the following lemma gives bounds of $|\nabla \rho|_q$, $|\nabla H|_q$ and $|\nabla^2 u|_q$.
 \begin{lemma}\label{s7}
 \begin{equation}\label{zhu54}
\begin{split}
&\|\big(\rho,H,P)(t)\|_{W^{1,r}}+|(\rho_t,H_t, P_t)(t)|_r\leq C,\quad \text{for} \quad 0\leq t\leq  T,
\end{split}
\end{equation}
where  $r\in[2,q]$, $C$ only depends on $C_0$, $\mu$, $\lambda$, $A$, $\gamma$, $\Omega$ and $T$ $(any \ T \in (0,\overline{T}])$.
\end{lemma}
\begin{proof}

In the following estimates we will use (from (\ref{diyi})-(\ref{diy2}) and (\ref{tvd9}))
 \begin{equation}\label{zhu55a}
\begin{split}
|\nabla^2 v|_q \leq& C(|\nabla \rho |_q+|\nabla H|_q),\\
|\nabla v|_{\infty}\leq& C\big(1+|\nabla v|_{BMO(\Omega)}\ln (e+|\nabla^2v|_q)\big)\\
\leq& C\big(1+(|\rho|_{L^2\cap L^\infty}+|H|_{L^2\cap L^\infty})\ln (e+|\nabla \rho|_q+|\nabla H|_q)\big)\\
\leq& C\big(1+\ln (e+|\nabla \rho|_q+|\nabla H|_q)\big).
\end{split}
\end{equation}

First,  applying $\nabla$ to  $(\ref{eq:1.2})_3$, multiplying the resulting equations by $q|\nabla \rho|^{q-2} \nabla \rho$, we have
\begin{equation}\label{zhu20cccc}
\begin{split}
&(|\nabla \rho|^q)_t+\text{div}(|\nabla \rho|^qu)+(q-1)|\nabla \rho|^q\text{div}u\\
=&-q |\nabla \rho|^{q-2}(\nabla \rho)^\top D( u) (\nabla \rho)-q \rho|\nabla \rho|^{q-2} \nabla \rho \cdot \nabla \text{div}u.
\end{split}
\end{equation}
Then integrating (\ref{zhu20cccc}) over $\Omega$, we immediately obtain
\begin{equation}\label{zhu200}
\begin{split}
\frac{d}{dt}|\nabla \rho|^q_q
\leq& C|D( u)|_\infty|\nabla \rho|^q_q+C|\nabla^2 u|_q|\nabla \rho|^{q-1}_q\\
\leq & C(|\nabla w|_\infty+|\nabla v|_\infty) |\nabla \rho|^q_q+C(|\nabla^2 w|_q+|\nabla^2 v|_q)|\nabla \rho|^{q-1}_q.
\end{split}
\end{equation}

Second, applying $\nabla$ to  $(\ref{eq:1.2})_1$, multiplying the resulting equations by $q\nabla H |\nabla H|^{q-2}$, we have
\begin{equation}\label{zhu20qs}
\begin{split}
&(|\nabla H|^2)_t-qA:\nabla H|\nabla H|^{q-2}+q B: \nabla H|\nabla H|^{q-2}+qC : \nabla H|\nabla H|^{q-2}=0,
\end{split}
\end{equation}
where
\begin{equation}\label{mou6ll}
\begin{split}
A=\nabla (H\cdot \nabla u)=&(\partial_j H \cdot \nabla u^i)_{(ij)}+(H\cdot \nabla \partial_j u^i)_{(ij)},\\
B=\nabla (u \cdot \nabla H)=&(\partial_j u\cdot \nabla H^i)_{(ij)}+ (u\cdot \nabla \partial_j H^i)_{(ij)},\\
C=\nabla(H \text{div}u)=&\nabla H \text{div}u+H \otimes \nabla \text{div}u.
\end{split}
\end{equation}

Then integrating (\ref{zhu20qs}) over $\Omega$, due to
\begin{equation}\label{btbt}
\begin{split}
&\int_{\Omega}  A: \nabla H |\nabla H|^{q-2}\text{dx}
\leq  C|\nabla u|_\infty |\nabla H|^q_q+C|H|_\infty|\nabla H|^{q-1}_q |u|_{D^{2,q}},\\
&\int_{\Omega}  B: \nabla H|\nabla H|^{q-2} \text{dx}\\
=&\int_{\Omega}   \sum_{i,j,k}  \partial_j u^k \partial_k H^i \partial_j H^i|\nabla H|^{q-2} \text{dx}+\int_{\Omega}  \sum_{i,j,k} u^k  (\partial_{kj} H^i \partial_j H^i)|\nabla H|^{q-2} \text{dx}\\
=&C|\nabla u|_\infty |\nabla H|^q_q+\frac{1}{2}\int_{\Omega}  \sum_{k=1}^3  u^k\Big(\partial_{k} |\nabla H|^2|\nabla H|^{q-2} \Big)\text{dx}\\
=&C|\nabla u|_\infty |\nabla H|^q_q+\frac{1}{q}\int_{\Omega}  \sum_{k=1}^3  u^k  \partial_{k} |\nabla H|^{q} \text{dx}
\leq C|\nabla u|_\infty |\nabla H|^q_q,\\
&\int_{\Omega}  C: \nabla H|\nabla H|^{q-2} \text{dx}
\leq C|\nabla u|_\infty |\nabla H|^q_q+C|H|_\infty|\nabla H|^{q-1}_q |u|_{D^{2,q}},
\end{split}
\end{equation}
 we quickly obtain the following estimate:
\begin{equation}\label{zhu21qs}
\begin{split}
\frac{d}{dt}|\nabla H|^q_q
\leq& C(|\nabla u|_\infty+1)|\nabla H|^q_q+C | u|_{D^{2,q}}|\nabla H|^{q-1}_q\\
\leq & C(|\nabla w|_\infty+|\nabla v|_\infty) |\nabla H|^q_q+C(|\nabla^2 w|_q+|\nabla^2 v|_q)|\nabla H|^{q-1}_q.
\end{split}
\end{equation}

Then from  (\ref{zhu55a}), (\ref{zhu200}),  (\ref{zhu21qs})  and Gronwall's inequality, we immediately have
\begin{equation}\label{zhu2kkh}
\begin{split}
&\frac{d}{dt}(|\nabla \rho|^q_q+|\nabla H|^q_q)\\
\leq & C(1+|\nabla w|_\infty+|\nabla v|_\infty) (|\nabla \rho|^q_q+|\nabla H|^q_q)+C|\nabla^2 w|_q(|\nabla \rho|^{q-1}_q+|\nabla H|^{q-1}_q)\\
\leq &C\big(1+\|\nabla w\|_{W^{1,q}}+\ln (e+|\nabla \rho|_q+|\nabla H|_q)\big)(|\nabla \rho|^q_q+|\nabla H|^q_q)\\
&+C|\nabla^2 w|_q(|\nabla \rho|^{q-1}_q+|\nabla H|^{q-1}_q).
\end{split}
\end{equation}
Via   (\ref{zhu2kkh}) and notations:
$$
f=e+|\nabla \rho|_q+|\nabla H|_q,\quad g=1+\|\nabla w\|_{W^{1,q}},
$$
 we quickly have
$$
f_t\leq Cgf +Cf\ln f+Cg,
$$
which, together with Lemma \ref{sk4zz} and Osgood Lemma,  implies that
$$
\ln f(t)\leq C,\quad \text{for} \quad 0\leq t\leq  T.
$$
Then we have obtained the desired estimate for $|\nabla \rho|_q+|\nabla H|_q$.  And the upper  bound of  $|\nabla \rho|_r+|\nabla H|_r$ can be deduced via the   H\"older's inequality.

Finally, the estimates for $\rho_t$ and $H_t$  can be obtained easily via the following relation:
\begin{equation}\label{ghtkk}
\begin{cases}
H_t=H \cdot \nabla u-u \cdot \nabla H-H\text{div}u,\\[6pt]
\rho_t=-u\cdot \nabla \rho-\rho\text{div} u,\ P_t=-u \cdot \nabla P-\gamma P \text{div}u,
\end{cases}
\end{equation}
and the estimates obtained in Lemmas  \ref{s2}-\ref{s7}.
 \end{proof}

According to the estimates obtained in Lemmas \ref{s2}-\ref{s7}, we deduce that
  \begin{lemma}\label{sk4nn}
\begin{equation*}
\begin{split}
|u(t)|_{D^2}+|\sqrt{\rho}u_t(t)|_{2}+\int_0^T\big(|u_t|^2_{D^1}+|u|^2_{D^{2,q}}\big) \text{d}t\leq C,\quad \text{for} \quad 0\leq t\leq  T,
\end{split}
\end{equation*}
 where $C$ only depends on $C_0$, $\mu$, $\lambda$, $A$, $\gamma$, $\Omega$ and $T$ $(any\  T\in (0,\overline{T}])$. \end{lemma}
\begin{proof}
Via the momentum equations $(\ref{eq:1.2})_4$, (\ref{diy2}) and Lemma \ref{tvd1},  we have
$$
|u|_{D^{2,l}}\leq (|w|_{D^{2,l}}+|v|_{D^{2,l}})\leq C(|w|_{D^{2,l}}+|\nabla P|_l+|\nabla H|_l),
$$
which, together with Lemma \ref{s7},  implies that
$$
|u(t)|_{D^2}+\int_0^T |u|^2_{D^{2,q}}\text{d}t\leq C,\quad \text{for} \quad 0\leq t\leq  T.
$$
According to Lemmas \ref{abs3} and \ref{s7},  for $r\in (3,7/2)$, we quickly have
\begin{equation*}
\begin{split}
|\sqrt{\rho}u_t|_2\leq C(|\sqrt{\rho}\dot{u}|_2+|\sqrt{\rho} u\cdot \nabla u|_2)\leq C(1+|\rho^{\frac{1}{r}}u|_r |\nabla u|_{\frac{2r}{r-2}}) \leq C.
\end{split}
\end{equation*}
Similarly, we have
$$\int_0^T |u_t|^2_{D^1} \text{d}t \leq C\int_0^T (|\dot{u}|^2_{D^1}+|u\cdot \nabla u|^2_{D^1})\text{d}t \leq C.$$
\end{proof}

In view of the estimates obtained in  Lemmas \ref{s2}-\ref{sk4nn}, we can extend the strong solutions of $(H,\rho,u,P)$ beyond $t\geq \overline{T}$. In truth, according to Lemmas \ref{s2}, \ref{sk4ss} and \ref{s7}-\ref{sk4nn},  the functions $(H,\rho,u,P)|_{t=\overline{T}} =\lim_{t\rightarrow \overline{T}}(H,\rho,u,P)$ satisfy the conditions imposed on the initial data (\ref{th78}) at the time $t=\overline{T}$. Furthermore,
\begin{equation}\label{th79mowei}
\begin{split}
(Lu+\nabla P- \text{rot} H\times H)|_{t=\overline{T}}=\lim_{t\rightarrow \overline{T}}(\rho u_t+\rho u\cdot \nabla u)=\sqrt{\rho} g|_{t=\overline{T}}
\end{split}
\end{equation}
with $g|_{t=\overline{T}}\in L^2$. Thus $(H,\rho,u,P)|_{t=\overline{T}}$ satisfy  $(\ref{th79})$ also. Then, we can take $(H,\rho,u,P)|_{t=\overline{T}}$ as the initial data and apply the local existence Theorem \ref{th5} to extend our local strong solution beyond $ \overline{T}$. This contradicts the assumption on $\overline{T}$.

{\bf Acknowledgement:} The research of  S. Zhu was supported in part
by National Natural Science Foundation of China under grant 11231006, Natural Science Foundation of Shanghai under grant 14ZR1423100 and  China Scholarship Council.

\end{document}